\documentclass[MSNbibl,number,citesort,dvips]{arxbj}
\usepackage{upgreek}
\usepackage{dcolumn}
\usepackage{graphicx}

\aid{0}
\volume{19}
\issue{3}
\pubyear{2013}
\firstpage{982}
\lastpage{1005}
\doi{10.3150/12-BEJ421} 

\makeatletter
\newcolumntype{d}[1]{D{.}{.}{#1}}
\newtheorem{them}[prop]{Theorem}
\newtheorem{lem}[prop]{Lemma}
\newtheorem{prop}{Proposition}[section]
\newtheorem{cor}[prop]{Corollary}
\newremark{rem}{Remark}

\renewcommand{\hat}{\widehat}
\newcommand{\bone}{{\boldsymbol{1}}}
\newcommand{\supp}{the supplemental article~\cite{veillette-Rosen-supp:2011}}
\makeatother

\begin{document}
\begin{frontmatter}

\title{Properties and numerical evaluation of the Rosenblatt distribution}
\runtitle{Properties and numerical evaluation of the Rosenblatt distribution}

\begin{aug}
\author{\fnms{Mark S.} \snm{Veillette}\thanksref{e1}\ead[label=e1,mark]{mveillet@bu.edu}} \and
\author{\fnms{Murad S.} \snm{Taqqu}\corref{}\thanksref{e2}\ead[label=e2,mark]{murad@bu.edu}}
\runauthor{M.S. Veillette and M.S. Taqqu} 
\address{Department of Mathematics and Statistics,
Boston University,
111 Cummington Street,
Boston, MA 02215, USA. \printead{e1,e2}}
\end{aug}

\received{\smonth{10} \syear{2010}}
\revised{\smonth{12} \syear{2011}}

%
\begin{abstract}
This paper studies various distributional properties of the Rosenblatt
distribution. We begin by describing a technique for computing the
cumulants. We then study the expansion of the Rosenblatt distribution
in terms of
shifted chi-squared distributions. We derive the coefficients of this
expansion and use these to obtain the L\'{e}vy--Khintchine formula and
derive asymptotic properties of the
L\'{e}vy measure. This allows us to compute the cumulants, moments,
coefficients in the chi-square expansion and the density and
cumulative distribution functions of the Rosenblatt distribution with
a high degree of precision. Tables are provided and software written
to implement the methods described here is freely available by request from
the authors.
\end{abstract}

%

\end{frontmatter}

\section{Introduction}

Typical limits of normalized of sums of long-range dependent stationary
series are Brownian motion, fractional Brownian motion or the
Rosenblatt process. Brownian motion and fractional Brownian motion are
Gaussian and can thus be readily tabulated. This is not the case for
the Rosenblatt distribution. The goal of this paper is to fill this
gap. The tables can be used to compute asymptotic confidence intervals
and to implement maximum likelihood methods.

The Rosenblatt distribution is the simplest
non-Gaussian distribution which arises in a
\textit{non-central limit theorem} involving long-range dependent
random variables~\cite{dobrushin:1979,taqqu:1975,taqqu:1979}. For an overview, see~\cite{taqqu:2011}. It also
appears in a statistical context as the asymptotic distribution of
certain estimators (e.g.,~\cite{tudor:2008}).

We shall begin by motivating the
Rosenblatt distribution using Rosenblatt's famous
counterexample found in~\cite{rosenblatt:1961}. Consider a stationary
Gaussian sequence $X_i$, $i = 1,2,\ldots$ which has a covariance
structure of the form $\mathbb{E} X_0 X_k \sim k^{-D}$ as $k
\rightarrow\infty$ with $0 < D <
1/2$. Using the transformation
%
\begin{equation}
Y_i = X_i^2 - 1,
\end{equation}
one can define a sequence of normalized sums
%
\begin{equation}\label{e:scalerv}
Z_{D}^{N} = \frac{\sigma(D)}{N^{1-D}} \sum_{i=1}^N Y_i.
\end{equation}
Here, $\sigma(D)$ is a normalizing constant and is given by
%
\begin{equation}\label{e:defofsigma}
\sigma(D) = \bigl[ \tfrac{1}{2} (1 - 2 D)(1- D)\bigr ]^{1/2}.
\end{equation}
The sequence $Z_D^{N}$ tends to a non-Gaussian limit $Z_D$ as $N
\rightarrow
\infty$ with mean $0$ and variance $1$. This limiting distribution was
named the Rosenblatt
distribution in~\cite{taqqu:1975}. The characteristic function of
$Z_D$ can be given as the following power series which is only
convergent near the origin:
%
\begin{equation} \label{e:charfuncsum}
\phi(\theta) = \exp\Biggl( \frac{1}{2} \sum_{k=2}^\infty(2 \mathrm{i} \theta
\sigma(D))^k \frac{c_k}{k} \Biggr),
\end{equation}
where
%
\begin{equation}\label{e:cumulantintegrals}
\hspace*{-5pt}c_k = \int_0^1 \mathrm{d}x_1 \int_0^1 \mathrm{d}x_2 \cdots\int_0^1 \mathrm{d}x_k |x_1 - x_2|^{-D}
|x_2 -x_3|^{-D} \cdots|x_{k-1} - x_k|^{-D} |x_k - x_1|^{-D}.
\end{equation}
By Cauchy--Schwarz,
%
\begin{equation}
\label{e:ck-comp}
\hspace*{-5pt}c_k \leq \biggl( \int_0^1 \int_0^1 |x_1 - x_2|^{-2D} \,\mathrm{d}x_1 \,\mathrm{d}x_2
\biggr)^{k/2}
= \biggl( \frac{1}{(1-2D)(1-D)}\biggr )^{k/2} = \biggl( \frac{1}{2 \sigma^2(D)} \biggr)^{k/2},
\end{equation}
ensuring that the series (\ref{e:charfuncsum}) converges around the origin.
[Since (\ref{e:ck-comp}) is an equality when $k=2$, $Z_D$ has variance $1$
in view of (\ref{e:charfuncsum}) and (\ref{e:defofsigma}).]

It is interesting to consider the extremes when $D\rightarrow0^+$ and $D
\rightarrow\frac{1}{2}^-$. When $D\rightarrow0^+$, notice that
$c_k \rightarrow1$ for all $k$, $\sigma(D) \rightarrow
\frac{1}{\sqrt{2}}$ and thus for
$\theta$ small enough, the characteristic function approaches
%
\begin{eqnarray}
\phi(\theta) &=& \exp\Biggl( \frac{1}{2} \sum_{k=2}^\infty\frac{(\sqrt
{2} \mathrm{i} \theta
)^k}{k} \Biggr) \nonumber\\
&=& \exp\biggl( \frac{1}{2} \bigl( \log\bigl(1 - \sqrt{2} \mathrm{i}
\theta\bigr) - \sqrt{2} \mathrm{i} \theta \bigr) \biggr) \\
&=& \biggl(\frac{1}{1
- \sqrt{2} \mathrm{i} \theta} \biggr)^{1/2} \mathrm{e}^{-\mathrm{i} \theta/\sqrt{2} },\nonumber
\end{eqnarray}
which is the characteristic function of $\frac{1}{\sqrt{2}} (\varepsilon
^2 - 1)$, where
$\varepsilon$ is $N(0,1)$. Hence when $D=0$, the Rosenblatt distribution
is simply a chi-squared distribution standardized to have mean 0 and
variance 1.

As $D\rightarrow\frac{1}{2}^-$, the limit is $N(0,1)$. This is not
surprising given that the scaling term in (\ref{e:scalerv}) approaches
$\sqrt{N}$, hence you would assume the usual central limit theorem to
hold and the limiting distribution to be Gaussian. This fact is not
obvious however from the characteristic function
(\ref{e:charfuncsum}). In Section~\ref{s:LKform7} of this work, we
derive an alternative form of
the characteristic function from which this Gaussian limit is easier to see.

The distribution $Z_D$ can be given in terms of a weighted sum of
chi-squared distributions,
%
\begin{equation}\label{e:chisum}
Z_D = \sum_{i=1}^\infty\lambda_n (\varepsilon_n^2 - 1),\qquad \varepsilon_n
\mbox{ i.i.d. }N(0,1),
\end{equation}
where the weights $\{\lambda_n\}_{n=1}^\infty$ are such that
%
\begin{equation}\label{e:ROSseries}
\sum_{n=1}^\infty\lambda_n^k = \sigma^k (D) c_k, \qquad k=2,3,\ldots.
\end{equation}
The series (\ref{e:chisum}) converges a.s. and in $L^2$
because
\begin{eqnarray*}
\sum_{n=1}^\infty\operatorname{Var} [\lambda_n (\varepsilon_n^2 - 1)]
&=& \mathbb{E}
[(\varepsilon_1^2 - 1)^2] \sum_{n=1}^\infty\lambda_n^2 = 2 \sum_{n=1}^\infty\lambda_n^2 < \infty.
\end{eqnarray*}
In fact, $\sum_{n=1}^\infty\lambda_n^2 = 1/2$ by (\ref
{e:cumulantintegrals}) and (\ref{e:defofsigma}).
The weights $\{\lambda_n\}$ are given as the eigenvalues of an
integral operator which we
will discuss in more detail in Section~\ref{s:evals}. Various integral
representations can be found in~\cite{tudor:2008}. Our main focus in
this work is on distributional properties of this
distribution, namely, cumulants, moments and obtaining a numerical
evaluation of the Rosenblatt distribution. A table for the cumulative
distribution function (CDF) of the Rosenblatt distribution is useful
for obtaining percentiles and confidence intervals.\looseness=-1

This paper is organized as follows: In Section~\ref{s:momcum}, we look
at the moments and cumulants of the Rosenblatt distribution, as well as
detailing
a method for computing them. In Section~\ref{s:evals}, we show that the
$\lambda_n$'s in the expansion (\ref{e:chisum}) are given by
the eigenvalues of an integral operator, and we give asymptotic
formulas for this sequence. In Section~\ref{s:LKform7}, we state the
characteristic function of
the Rosenblatt distribution in L\'{e}vy--Khintchine form and use it to
derive further properties of the Rosenblatt distribution.
In Section~\ref{s:compute}, we compute the moments,
cumulants and the $\lambda_n$'s and in Section~\ref{s:computeCDF}, the
quantiles of the Rosenblatt distribution are computed for various $D$ values.
For more details and software, see \supp.
%
\section{Cumulants and moments of the Rosenblatt distribution}\label{s:momcum}
It follows from the expansion of the characteristic function (\ref
{e:charfuncsum}) around
$\theta= 0$ that the cumulants $\kappa_k$
of the Rosenblatt distribution are given by $\kappa_1 = 0$ and
%
\begin{equation}\label{e:cumulantformula}
\kappa_k = 2^{k-1} (k-1)! (\sigma(D))^k c_k,\vadjust{\goodbreak}
\end{equation}
where the $c_k$ are given by the multiple integrals
(\ref{e:cumulantintegrals}). Each moment $\mu_n$, $n \geq1$ can then be
expressed as a polynomial using the cumulants $\kappa_k$,
$k=1,2,\ldots,n$. These are the \textit{complete Bell Polynomials}
%
\begin{equation}\label{e:momentsBell}
\mu_n = B_n(0,\kappa_2,\ldots,\kappa_n)
\end{equation}
as noted, for example, in~\cite{taqqu:2010,pitman:2002}. They
can also be computed recursively~\cite{smith:1995}.

Thus, in order to compute any moment or cumulant, it is necessary to
compute the multiple integrals $c_k$. The first two can be computed
directly,
%
\begin{eqnarray}\label{e:c3}
c_2 &=& \int_0^1 \int_0^1 |x_1 - x_2|^{-2D} \,\mathrm{d}x_1 \,\mathrm{d}x_2 = 2 \int_0^1
\int_0^{x_2} (x_1 - x_2)^{-2D} \,\mathrm{d}x_1 \,\mathrm{d}x_2 =
\frac{1}{(1-2D)(1-D)}, \nonumber\\
c_3 &=& \int_0^1 \int_0^1 \int_0^1 |x_1 - x_2|^{-D} |x_2 - x_3|^{-D}
|x_3 - x_1|^{-D} \,\mathrm{d}x_1 \,\mathrm{d}x_2 \,\mathrm{d}x_3 \nonumber\\
&=& 3 \int_0^1 x_3^{3D} \int_0^{x_3} \int_0^{x_3} \biggl|\frac{x_1}{x_3} -
\frac{x_2}{x_3} \biggl|^{-D} \biggl(1 - \frac{x_2}{x_3} \biggr)^{-D} \biggl(1 -
\frac{x_1}{x_3} \biggr)^{-D} \,\mathrm{d}x_1 \,\mathrm{d}x_2 \,\mathrm{d}x_3 \nonumber\\
&=& 3 \biggl( \int_0^1 x_3^{-3D+2} \,\mathrm{d}x_3 \biggr) \biggl( \int_0^{1} \int_0^{1} |u_1 -
u_2|^{-D} (1
- u_2)^{-D} (1- u_1)^{-D} \,\mathrm{d}u_1 \,\mathrm{d}u_2\biggr ) \\
&=& \frac{2}{1 - D} \int_0^1 w_2^{-D} \int_0^{w_2} (w_2 -
w_1)^{-D} w_1^{-D} \,\mathrm{d}w_1 \,\mathrm{d}w_2 \nonumber\\
&=& \frac{2}{1-D} \biggl( \int_0^1 w_2^{-3D + 1} \,\mathrm{d}w_2 \biggr) \biggl(
\int_0^1 v^{-D} (1-v)^{-D} \,\mathrm{d}v \biggr) \nonumber\\
&=& \frac{2}{(1-D)(2 - 3D)} \beta(1-D,1-D),\nonumber
\end{eqnarray}
where $\beta(a,b) = \int_0^1 v^{a - 1}(1-v)^{b-1}\,\mathrm{d}v$ is the beta
function and we made the following changes of variables above: $u_1 =
x_1/x_3, u_2 =
x_2/x_3$, $w_1 = 1 - u_1,w_2 = 1-u_2$ and $v = w_1/w_2$.

For $k \geq4$, a closed form expression for $c_k$ could not be
found, which means they must be computed numerically. Computing
the multiple integrals directly is intractable due to the increasing
number of singularities in the integrand. It is for this reason that we
now develop a more sophisticated method for computing $c_k$.

Let $L^2(0,1)$ denote the Hilbert space of all real-valued measurable
functions $f(x), 0 < x < 1$, such that $\Vert f\Vert_2 \equiv(\int_0^1 f(x)^2
\,\mathrm{d}x)^{1/2} < \infty$, together with the usual inner product $(f,g)
\equiv
\int_0^1 f(x) g(x) \,\mathrm{d}x$. For $ 0 < D < 1/2$, define the integral
operator $\mathcal{K}_D
\dvtx L^2(0,1) \rightarrow L^2(0,1)$ as
%
\begin{equation}\label{e:Koperator}
(\mathcal{K}_D f) ( x) = \int_0^1 |x - u|^{-D} f(u) \,\mathrm{d}u.
\end{equation}
Finally, define the sequence of functions $G_{k,D}\in L^2(0,1)$, $ k
\geq1$,
recursively as follows:
%
\begin{equation}\label{e:defofGkD}
G_{1,D}(x) = \frac{(1-x)^{-D}}{\sqrt{1-D}};\qquad
G_{k,D}(x) = (\mathcal{K}_D G_{k-1,D})(x), \qquad k \geq2.
\end{equation}
Then, we have the following alternative way to express $c_k$.
\begin{prop}\label{p:computeint}
Let $\mu$ and $\nu$ be any two positive integers such that $\mu+ \nu
= k$. Then
%
\begin{equation}\label{e:computeint}
c_k = (G_{\mu,D},G_{\nu,D}).
\end{equation}
\end{prop}
\begin{pf}
Let $\mu,\nu$ be as stated. Then, using the circular symmetry of the
integrand in (\ref{e:cumulantintegrals}), if we take $x_k$ as the
largest of the $x_i$, $i=1,2,\ldots$\,, and then factor an $x_k$ out of
all the terms, we can rewrite $c_k$ as
%
\begin{eqnarray}
c_k &=& k \int_0^1 x_k^{-k D} \int_{(0,x_k)^{k-1}} \biggl( 1-
\frac{x_1}{x_k} \biggr)^{-D} \biggl| \frac{x_1}{x_k} - \frac{x_2}{x_k}
\biggr|^{-D} \ldots\nonumber
\\[-8pt]
\\[-8pt]
&&\hphantom{k \int_0^1 x_k^{-k D} \int_{(0,x_k)^{k-1}}}{} \times \biggl| \frac{x_{k-2}}{x_k} - \frac{x_{k-1}}{x_k}
\biggr|^{-D} \biggl( 1 - \frac{x_{k-1}}{x_k} \biggr) ^{-D} \,\mathrm{d}x_1 \cdots \,\mathrm{d}x_k.\nonumber
\end{eqnarray}
With the change of variables $u_i = x_i/x_k$, $i = 1,2,\ldots,k-1$, one
of the $k$ integrals can be separated out, and we obtain
\begin{eqnarray*}
 c_k &=& k \biggl(\int_0^1 x_k^{-k D + (k-1)} \,\mathrm{d}x_k \biggr) \\
&&{} \times \biggl(
\int_{(0,1)^{k-1}} (1-u_1)^{-D} |u_1 - u_2|^{-D} \cdots
|u_{k-1} -
u_{k-2}|^{-D} (1- u_{k-1})^{-D} \,\mathrm{d}u_1 \cdots \,\mathrm{d}u_{k-1} \biggr)
\\
&=& \frac{1}{1 - D} \biggl(
\int_{(0,1)^{k-1}} (1-u_1)^{-D} |u_1 - u_2|^{-D} \cdots\\
&&\hphantom{\frac{1}{1 - D} \biggl(
\int_{(0,1)^{k-1}}}{}\times|u_{k-1} -
u_{k-2}|^{-D} (1- u_{k-1})^{-D} \,\mathrm{d}u_1 \cdots \,\mathrm{d}u_{k-1} \biggr) \\
&=& \int_{(0,1)^{k-1}} G_{1,D}(u_1) \underbrace{|u_1 - u_2|^{-D}
\cdots|u_{k-1} -
u_{k-2}|^{-D} }_{k - 2 \ \mathrm{terms}}G_{1,D}(u_{k-1}) \,\mathrm{d}u_1 \cdots \,\mathrm{d}u_{k-1} \\
&=& \int_{(0,1)^{k-3}} \biggl[ \int_0^1 G_{1,D}(u_1) |u_1 - u_2|^{-D}
\,\mathrm{d}u1 \biggr] [ |u_3 - u_2|^{-D} \cdots|u_{k-2} - u_{k-3}|^{-D}
]\\
&&\hphantom{\int_{(0,1)^{k-3}}}{} \times
\biggl[ \int_0^1 G_{1,D}(u_{k-1}) |u_{k-1} - u_{k-2}|^{-D}
\,\mathrm{d}u_{k-1}\biggr ] \,\mathrm{d}u_3 \cdots \,\mathrm{d}u_{k-2} \\
& =& \int_{(0,1)^{k-3}} G_{2,D}(u_2) \underbrace{|u_2 - u_3|^{-D}
\cdots|u_{k-2} -
u_{k-3}|^{-D}}_{k - 4 \ \mathrm{terms}} G_{2,D}(u_{k-2}) \,\mathrm{d}u_2 \cdots \,\mathrm{d}u_{k-2} \\
& \vdots&\\
&=& \int_0^1 G_{\mu,D}(u_{\mu}) G_{\nu,D} (u_{k - \nu}) \,\mathrm{d}u_{\mu} =
(G_{\mu,D},G_{\nu,D}).
\end{eqnarray*}
This finishes the proof.
\end{pf}

\begin{rem*} To minimize the number of integrals one needs to
compute, it makes
sense to choose $\mu= \nu= \frac{k}{2}$ if $k$ is even, and
$\mu= \frac{k+1}{2}$ and $\nu= \frac{k-1}{2}$ if $k$ is odd.
Proposition~\ref{p:computeint} thus reduces the problem of computing a
$k$-dimensional integral into computing $\lceil\frac{k}{2} \rceil+
1$ one-dimensional integrals.
\end{rem*}

$G_{2,D}$ can be given in terms of the beta function and the Gauss
hypergeometric function ${}_2F_1(a,b;c;x)$, which has the following integral
representation,
%
\begin{eqnarray}\label{e:hypgeoproperty}
&&{}_2F_1(a,b;c;x) = \frac{\Gamma(c)}{\Gamma(b)\Gamma(c - b)} \int_0^1
v^{b-1} (1-v)^{c - b - 1} (1- x v)^{-a} \,\mathrm{d}v, \nonumber
\\[-8pt]
\\[-8pt]
&&\quad x < 1, c>b>0.\nonumber
\end{eqnarray}
Indeed,
%
\begin{eqnarray}\label{e:G2closed}
G_{2,D}(x) &=& \frac{1}{\sqrt{1-D}} \int_0^1 (1 - u)^{-D} |x - u|^{-D}
\,\mathrm{d}u \nonumber\\
&=& \frac{1}{\sqrt{1-D}} \biggl[ \int_0^x (1-u)^{-D} (x - u)^{-D} \,\mathrm{d}u +
\int_x^1 (1-u)^{-D} (u-x)^{-D} \,\mathrm{d}u \biggr] \nonumber\\
&=& \frac{1}{\sqrt{1-D}} \biggl[ x^{1 - D} \int_0^1 (1 - x v)^{-D} (
1- v)^{-D} \,\mathrm{d}v \nonumber
\\[-8pt]
\\[-8pt]
&&\hphantom{\frac{1}{\sqrt{1-D}} \biggl[}{} + (1-x)^{1-2D} \int_0^1 w^{-D} (1-w)^{-D} \,\mathrm{d}w \biggr] \nonumber
\\
&=& \frac{1}{\sqrt{1-D}} \biggl[ \frac{x^{1-D}}{1-D} {}_2F_1(D,1,2-D,x)
+ (1-x)^{1 - 2D} \beta(1-D,1-D) \biggr] \nonumber\\
&=& \frac{x^{1-D}}{(1-D)^{3/2}} {}_2F_1(D,1,2-D,x)
+ \frac{(1-x)^{1 - 2D} \beta(1-D,1-D)}{\sqrt{1-D}},\nonumber
\end{eqnarray}
where, in the third equality we used the change of variables $v = u/x$
and $w = (1 - u)/(1-x)$, and in the fourth, we used
(\ref{e:hypgeoproperty}). The function ${}_2F_1(a,b;c;x)$ is bounded
for $x \in(0,1)$ as long as $c > a+b$
(\cite{oldham:2009}, Section 60:7), which is true
in this case. This implies that unlike $G_{1,D}$, $G_{2,D}$ is
a bounded function on $(0,1)$, since $0 < D < \frac{1}{2}$.

In Section C
of \supp, we outline a technique for
computing the $c_k$ numerically based on Proposition
\ref{p:computeint}, and tabulate the first 8 cumulants and moments of
the Rosenblatt distribution for various values of $D$.
%
\section{Eigenvalue expansion of the Rosenblatt distribution}\label{s:evals}
In this section, we focus on the expansion of the Rosenblatt
distribution in terms of shifted chi-squared distributions
(\ref{e:chisum}). The sequence $\{\lambda_i\}_{i=1}^\infty$ can be
thought of in two ways.

One way is to start with the integrals $\{c_k, k \geq2\}$ defined in
(\ref{e:cumulantintegrals}) and to view $\{\lambda_i\}_{i=1}^\infty$
as a non-increasing sequence related to these $\{c_k, k \geq2\}$
through formula (\ref{e:ROSseries}),
see~\cite{taqqu:1975}.
While easier to state, this perspective sheds little
light on the $\lambda_i$'s since the $c_k$'s are so complicated.

The
second way to characterize the sequence $\{\lambda_i\}_{i=1}^\infty$
is more useful in our case, and stems from Proposition 2 in
\cite{dobrushin:1979}. To recall this proposition, let $X$ be
defined through the Wiener--It\^{o} integral
%
\begin{equation}\label{e:WIintegral}
X = \int^{\prime\prime}_{\mathbb{R}^2} H(x,y) Z_G(\mathrm{d}x) Z_G(\mathrm{d}y),
\end{equation}
where $Z_G$ is a complex-valued random measure with control measure $G$
such that for all Borel sets $A \in
\mathbb{R}$, $Z_G(A) = Z_G(-A)$ and $G(A) = G(-A)$, and the kernel
$H(x,y)$ is a complex-valued measurable function such that
$H(x,y) = H(y,x) = \overline{H(-x,-y)}$ for all $x,y \in\mathbb{R}$,
and $\int_{\mathbb{R}^2} |H(x,y)|^2 G(\mathrm{d}x)G(\mathrm{d}y) < \infty$. The
double prime in the integral means to exclude the diagonal $\{x =
\pm y\}$. For background on such integrals, see~\cite{taqqu:2010} or
\cite{major:1981}.

Let $L_G^2(\mathbb{R})$ denote the space of complex valued functions
$h(x)$, $x \in\mathbb{R}$ such that for all $x \in\mathbb{R}$, $h(x)
= \overline{h(-x)}$ and $\int|h(x)|^2 G(\mathrm{d}x) < \infty$. For random
variables $X$ defined as in~(\ref{e:WIintegral}), Dobrushin
and Major showed that $X$ has an expansion
\[
X = \sum_{n=1}^\infty\eta_n (\varepsilon_n^2 - 1),
\]
where the sequence $\eta_n$ corresponds to the eigenvalues of the
integral operator $A\dvtx L_G^2(\mathbb{R}) \rightarrow
L_G^2(\mathbb{R}^2)$ defined as
%
\begin{equation}
(A h)(x) = \int_{-\infty}^\infty H(x,-y) h(y) G(\mathrm{d}y).
\end{equation}

In~\cite{taqqu:1979}, it is shown that the
Rosenblatt distribution has the following representation as a
Wiener--It\^{o} integral:
%
\begin{equation}
Z_D = a(D) \int_{\mathbb{R}^2}^{\prime\prime} \frac{ \mathrm{e}^{\mathrm{i} (x + y)}
- 1}{\mathrm{i} (x + y)}
Z_{G_D}(\mathrm{d}x) Z_{G_D}(\mathrm{d}y),
\end{equation}
where the measure $G_D$ is absolutely continuous and is given by
%
\begin{equation}\label{e:measureG}
G_D(\mathrm{d}x) = |x|^{D-1} \,\mathrm{d}x,\qquad x \in\mathbb{R}.
\end{equation}
The constant
%
\begin{equation}\label{e:defofRosa}
a(D) = \frac{\sigma(D)}{2 \Gamma(D) \sin((1-D)\uppi/2)}
,
\end{equation}
where $\sigma(D)$ is given in (\ref{e:defofsigma}), ensures a
variance of 1.

Thus, Dobrushin and Major's result implies that the sequence
$\lambda_n$ we seek in (\ref{e:chisum}) is given as the eigenvalues of
the operator $A_D\dvtx L^2_{G_D}(\mathbb{R}) \rightarrow L^2_{G_D}(\mathbb
{R})$ defined as
%
\begin{equation}\label{e:oldoperator}
(A_D h)(x) = a(D) \int_{-\infty}^\infty\frac{\mathrm{e}^{\mathrm{i} (x - y)} - 1}{\mathrm{i}
(x-y) }
h(y) G_D(\mathrm{d}y).
\end{equation}

We shall now reexpress the eigenvalue problem associated to
(\ref{e:oldoperator}) in a much simpler form so that we may both give
analytical results about the $\lambda_n$'s and develop a method to
compute them. This is done in the following proposition.
\begin{prop}\label{p:sameevals}
The operators $A_D\dvtx L^2_{G_D}(\mathbb{R}) \rightarrow L^2_{G_D}(\mathbb
{R})$ defined in (\ref{e:oldoperator})
and $\sigma(D) \mathcal{K}_D\dvtx\break L^2(0, 1) \rightarrow L^2(0,1)$ defined
in (\ref{e:Koperator}) have the same eigenvalues.
\end{prop}
\begin{pf}
Let $\hat{g}(z)= (\mathcal{F} g)(z) = \int_{\mathbb{R}} g(y) \mathrm{e}^{\mathrm{i} y z}
\,\mathrm{d}x$ and $\check{g}(z) = (\mathcal{F}^{-1} g)(z) = \frac{1}{2 \uppi}
\int_{\mathbb{R}} \mathrm{e}^{- \mathrm{i} x z} g(y) \,\mathrm{d}y$ denote the Fourier transform
and inverse Fourier transfrom. Recall that $\mathcal{F}$ and $\mathcal
{ F}^{-1}$ are defined on $L^1(\mathbb{R})$ and $L^2(\mathbb{R})$ and
can also be extended to generalized functions (\cite{Zemanian:1987},
Chapter 7).
%

Let $(\lambda,h)$ be an
eigenpair of the operator $A_D$. This implies that $h \in
L_{G_D}^2(\mathbb{R})$ and thus $\int
h(y)^2 |y|^{D-1} \,\mathrm{d}y < \infty$ by (\ref{e:measureG}).
Taking inverse Fourier transforms of both sides in $\lambda h = A_D h$,
we obtain
%
\begin{equation}\label{e:H1convH2}
\lambda\check{h} = \mathcal{F}^{-1} (A_D h) = a(D) \mathcal{F}^{-1}
(H_1 * H_2),
\end{equation}
where
\[
H_1(y) = \bigl(\exp(\mathrm{i} y) - 1\bigr)/(\mathrm{i}y)
\]
and
\[
H_2(y) = |y|^{D-1} h(y).
\]
We want to apply the convolution theorem to compute the inverse Fourier
transform in~(\ref{e:H1convH2}), however some care must be taken,
because, while $H_1 \in L^1(\mathbb{R}) \cap L^2(\mathbb{R})$,
$H_2$ is not necessarily in either $L^1(\mathbb{R})$ or
$L^2(\mathbb{R})$. However, this difficulty can be avoided by writing
$H_2$ as a sum:
\[
H_2(y) = |y|^{D-1} h(y) \mathbf{1}_{[-1,1]}(y) + |y|^{D-1} h(y)
\mathbf{1}_{(1,\infty)}(|y|) := H_2^-(y) + H_2^+(y).
\]
Since
$
\int h(y)^2 |y|^{D-1} \,\mathrm{d}y < \infty
$,
we have
\[
H_2^-(y) = |y|^{D-1}
h(y) \mathbf{1}_{[-1,1]}(y) \in L^1(\mathbb{R})
\]
and
\[
H_2^+(y) =|y|^{D-1} h(y)
\mathbf{1}_{(1,\infty)}(|y|) \in L^2(\mathbb{R}).
\]
By
linearity of the convolution and Fourier transform, we can apply the
convolution theorem for $L^1$ functions (\cite{stade:2005}, Proposition
6.2.1) and the convolution theorem for $L^2$ functions (\cite{stade:2005},
Proposition 6.5.3), we have
\begin{eqnarray*}
\lambda\check{h} &=& a(D) \mathcal{F}^{-1} (H_1 * H_2) = a(D) [
\mathcal{F}^{-1} (H_1 * H_2^-) + \mathcal{F}^{-1} (H_1 * H_2^+)
] \\
&=& 2 \uppi a(D) \mathbf{1}_{(0,1)} (\check{H}_2^- + \check{H}_2^+) = 2
\uppi a(D) \mathbf{1}_{(0,1)} \check{H}_2,
\end{eqnarray*}
where we have used the fact that $\check{H_1} =
\mathbf{1}_{(0,1)}$ and we have picked up an extra factor of $2 \uppi$
since we are applying the inverse Fourier transform $\mathcal{F}^{-1}$
to the convolution. This implies that for any eigenfunction $h$ of
$A_D$, the support of
$\check{h}$ is contained in $(0,1)$. Viewing $H_2$ as the product of
$|y|^{D-1}$ and $h(y)$, we again apply the convolution theorem.
Again, care must be taken as $|y|^{D-1}$ is not integrable or square
integrable, but if we
view $h$ and $|y|^{D-1}$ as generalized functions, the convolution
theorem still applies since $\check{h}$ has compact support
(\cite{Zemanian:1987}, Theorem 7.9-1). Thus,
\[
h(y) |y|^{D-1} = \mathcal{F} \bigl( \check{h} * \mathcal{F}^{-1} (|y|^{D-1})
\bigr).
\]
Using this, we have
\begin{eqnarray*}
2 \uppi a(D) \mathbf{1}_{(0,1)} \check{H}_2 &=& 2 \uppi a(D)
\mathbf{1}_{(0,1)} \mathcal{F}^{-1} [ h(y) |y|^{D-1} ] \\
&=& 2 \uppi a(D)
\mathbf{1}_{(0,1)} \mathcal{F}^{-1} \bigl[ \mathcal{F}\bigl ( \check{h} *
\mathcal{ F}^{-1}(|y|^{D-1}) \bigr) \bigr] \\
&=& 2 \uppi a(D)
\mathbf{1}_{(0,1)} \bigl( \check{h} * \mathcal{ F}^{-1}(|y|^{D-1}) \bigr).
\end{eqnarray*}
The inverse Fourier transform of $|y|^{D-1}$ is given by (\cite
{gradshteyn:2007}, page 1119),
%
\begin{equation}
(\mathcal{F}^{-1} |y|^{D-1})(z) = \frac{1}{2 \uppi} \biggl[2 \Gamma(D) \sin
\biggl( \frac{(1-D) \uppi}{2}\biggr )
|z|^{-D} \biggr]= \frac{\sigma(D)}{2 \uppi a(D)} |z|^{-D}
\end{equation}
from (\ref{e:defofRosa}). Thus, if $(h,\lambda)$ is an eigenpair for
$A_D$, then
%
\begin{eqnarray}\label{e:lamhhat1}
\lambda\check{h}(z) &=& (\mathcal{F}^{-1} A_D h) (z) \nonumber\\
&=& \sigma(D) \bone_{[0,1]}(z)
\int_{-\infty}^\infty|z - y|^{-D} \check{h}(y) \,\mathrm{d}y \\
&=& \sigma(D) \bone_{[0,1]}(z)
\int_{0}^1 |z - y|^{-D} \check{h}(y) \,\mathrm{d}y,\nonumber
\end{eqnarray}
where we have again used the fact that $\check{h}$ is supported on
$(0,1)$. Thus, if $(\lambda,h)$ is an eigenpair for $A_D$, then
$(\lambda,\check{h}|_{(0,1)})$ is an eigenpair for $\sigma(D)
\mathcal{ K}_D$. Reversing this argument shows that there is a one-to-one
correspondence between eigenpairs of $A_D$ and $\sigma_D \mathcal{K}$,
which preserves the eigenvalues, hence these operators have the same
eigenvalues. This completes the proof.
\end{pf}

The eigenvalues $\lambda_n$ of $\mathcal{K}_D$ are not known exactly,
but their
asymptotic behavior is well understood as $n \rightarrow\infty$, see
for instance
\cite{dostanic:1998,birman:1977,kac:1955} or
\cite{rosenblatt:1963}. In particular, we have the following result
regarding the
asymptotic behavior of $\lambda_n$ in (\ref{e:chisum}).
\begin{them}\label{t:evals}
Let $Z_D$ denote the Rosenblatt distribution given by
%
\begin{equation}\label{e:sumchisinEVALTHM}
Z_D = \sum_{n=1}^\infty\lambda_n(D) (\varepsilon_n^2 - 1),\qquad \varepsilon_n
\mbox{ i.i.d. } N(0,1),
\end{equation}
where the non-increasing sequence $\{\lambda_n(D)\}$ is given by the
eigenvalues of the integral operator $\sigma(D) \mathcal{K}_D\dvtx L^2(0,1)
\rightarrow L^2(0,1)$.
The following asymptotic formula holds for any $0 < r < 1$:
%
\begin{equation}\label{e:asylambdas}
\lambda_n(D) = C(D) n^{D-1} \biggl( 1 + \mathrm{o} \biggl(\frac{1}{n^r} \biggr) \biggr),
\end{equation}
where,
%
\begin{equation}\label{e:defofC}
C(D) = \frac{2}{\uppi^{1-D}} \sigma(D) \Gamma(1-D)
\sin\biggl(\frac{\uppi D}{ 2} \biggr) .
\end{equation}
Moreover, the series
\[
\sum_{n=1}^\infty \bigl(\lambda_n(D) - C(D) n^{D-1}\bigr)
\]
converges and equals
%
\begin{equation}\label{e:diffintailasy}
\sum_{n=1}^\infty \bigl(\lambda_n(D) - C(D) n^{D-1}\bigr) = -2^{1-D} \sigma(D)
\zeta(D),
\end{equation}
where $\zeta$ denotes the Riemann zeta function.
\end{them}
\begin{pf}
Proposition 2 in~\cite{dobrushin:1979} and Proposition
\ref{p:sameevals} verify the first claim, namely that
(\ref{e:sumchisinEVALTHM}) holds and that the $\lambda_n$'s are the
eigenvalues of $\sigma(D) \mathcal{K}_D$. Theorem 1 in
\cite{dostanic:1998} describes the eigenvalues of the operator
$\int_{-1}^1 |x - u|^{-D} f(u) \,\mathrm{d}y$. Thus, (\ref{e:asylambdas})--(\ref{e:diffintailasy}) follow
immediately after noting that
\[
\int_{-1}^1 |z - u|^{-D}
g(u) \,\mathrm{d}u= 2^{1-D} \int_{0}^1 |x - y|^{-D} f(y) \,\mathrm{d}y,
\]
where $z = 2x - 1$, $u = 2y - 1$ and $g(u) =
f (\frac{u+1}{2} )$.
\end{pf}

Since Theorem~\ref{t:evals} only gives the asymptotic behavior of the
$\lambda_n$'s, we need to approximate those for which the
asymptotic formula is not applicable. We present a method for
approximating these eigenvalues in Section~\ref{s:compute} (see also
Section D
of \supp).
%
\section{L\'{e}vy--Khintchine representation of the Rosenblatt
distribution}\label{s:LKform7}
Recall that a distribution $X$ is infinitely divisible if for any
integer $n
\geq1$, there exits $X_i^{(n)}$, $i=1,2,\ldots,n$, i.i.d. such that
$X \stackrel{d}{=} X_1^{(n)} + X_2^{(n)} + \cdots + X_n^{(n)}$. The
characteristic function of an infinitely divisible
distribution $X$ with $\mathbb{E} X^2 < \infty$ can always be written
in the following form:
\[
\phi(\theta) = \mathbb{E} \mathrm{e}^{\mathrm{i} \theta X} = \exp\biggl( \mathrm{i} a \theta-
\frac{1}{2} b^2 \theta^2 +
\int_{-\infty}^\infty(\mathrm{e}^{\mathrm{i} \theta u} - 1 - \mathrm{i} u \theta) \nu(\mathrm{d}u)
\biggr),
\]
where $a \in\mathbb{R}$, $b > 0$ and $\nu$ is a positive measure on
$\mathbb{R} \setminus\{0\}$ with
the property that $\int\min(u^2,1)\* \nu(\mathrm{d}u) < \infty$. This is known
as the L\'{e}vy Khintchine representation of $X$. Since the chi-square
distribution is infinitely divisible, it is not
surprising in light of (\ref{e:chisum}) that the Rosenblatt
distribution is also infinitely divisible. In this section, we will
make this assertion rigorous, and
give the L\'{e}vy Khintchine representation of the Rosenblatt
distribution.\looseness=-1

Before stating the result, we require first a lemma. Given any
positive, increasing
sequence $\mathbf{c} = \{c_n\}_{n=1}^\infty$ such that $\sum1/c_n^2 <
\infty$, define the
function $G_{\mathbf{c}}(x)$ for $0 < x < 1$ as
%
\begin{equation}\label{e:defofG}
G_{\mathbf{c}} (x) = \sum_{n=1}^\infty x^{c_n}.
\end{equation}
Since $\mathbf{c}^{-1} = \{c_1^{-1},c_2^{-1},\ldots\} \in\ell^2$,
we have $c_n \rightarrow\infty$ and thus this
series converges for all $x \in(0,1)$ since $\log x < 0 $ and hence
for $n$ large enough,
one has
$
x^{c_n} = \mathrm{e}^{c_n \log(x)} \leq c_n^{-2}.
$
Notice that $G(0) = 0$, $G_{\mathbf{c}}(x) \rightarrow\infty$ as $x
\rightarrow
1$ and is a continuous function for all $x \in(0,1)$.

We now state a lemma regarding the asymptotic behavior of
$G_{\mathbf{c}}$ as $x \rightarrow0$ and $x \rightarrow1$. In the
following, we will say $a_n
\sim b_n$ if $a_n/b_n \rightarrow1$, and $a_n \lesssim b_n$ if
$\limsup a_n/b_n \leq1$.
\begin{lem}\label{l:gengeosum}
Suppose $\mathbf{c}$ is a positive strictly increasing sequence such
that $c_n \sim\beta n^{\alpha}$
as $n \rightarrow\infty$ for some $1/2 < \alpha< 1$ and constant
$\beta
> 0$. Then,
%
\begin{eqnarray}
\label{e:asyG0}G_{\mathbf{c}}(x) &\sim& x^{c_1}, \qquad\mbox{as } x \rightarrow0,\\
\label{e:asyGinf}G_{\mathbf{c}}(x) &\sim&\frac{1}{\alpha\beta^{1/\alpha}}
\Gamma\biggl( \frac{1}{\alpha} \biggr)
(1-x)^{-1/\alpha}, \qquad\mbox{as } x \rightarrow1.
\end{eqnarray}
\end{lem}
\begin{pf}
As $x \rightarrow0$, we have
%
\begin{equation}
\frac{G_{\mathbf{c}}(x)}{x^{c_1}} = 1 + \sum_{n=2}^\infty x^{c_n -
c_1} \rightarrow1,
\end{equation}
since $c_n > c_1$, so that the sum on the right hand side tends to 0.
This confirms (\ref{e:asyG0}).\vadjust{\goodbreak}

For the second assertion, let $0 < \varepsilon<1 $ and let $\beta' =
(1-\varepsilon)\beta$ and $\beta'' = (1+\varepsilon) \beta$. By
assumption, there exists $M$ large enough
such that for $n \geq M$,
%
\begin{equation}
\beta' n^\alpha< c_n < \beta'' n^\alpha.
\end{equation}
Let $G_{\mathbf{c}}^{(M)}(x)$ be the tail of the series
$G_{\mathbf{c}}$:
\[
G_{\mathbf{c}}^{(M)}(x) = \sum_{n=M}^\infty x^{c_n}.
\]
Now, observe,
%
\begin{equation}\label{e:inqstring}
\int_M^\infty x^{(1+\varepsilon) \beta(y+1)^\alpha} \,\mathrm{d}y \leq\sum
_{n=M}^\infty
x^{\beta'' n^\alpha} \leq G_{\mathbf{c}}^{(M)}(x) \leq\sum
_{n=M}^\infty
x^{\beta' n^\alpha} \leq\int_M^\infty x^{(1-\varepsilon) \beta
y^\alpha}
\,\mathrm{d}y.
\end{equation}
The integrals on the far ends of this inequality can be computed in
terms of the upper incomplete gamma function $\Gamma(a,b) =
\int_b^\infty z^{a-1} \mathrm{e}^{-z} \,\mathrm{d}z$, since for the right most integral,
%
\begin{eqnarray}
\int_M^\infty x^{ \beta' y^\alpha} \,\mathrm{d}y &=& \int_M^\infty
\mathrm{e}^{ \beta' y^\alpha\log(x) } \,\mathrm{d}y \\
&=&\frac{1}{\alpha} \frac{1}{ (-\beta' \log(x) )^{1/\alpha}
} \int_{-\beta' \log(x)
M^{\alpha} }^\infty u^{1/\alpha- 1} \mathrm{e}^{-u} \,\mathrm{d}u \\
\label{e:rightside}&=& \frac{(-\log(x))^{ -1/\alpha}}{\alpha\beta'^{1/\alpha} }
\Gamma\biggl( \frac{1}{\alpha} , -\beta' \log(x) M^\alpha \biggr).
\end{eqnarray}
And similarly, for the left most side of (\ref{e:inqstring}),
%
\begin{equation}\label{e:leftside}
\int_M^\infty x^{ \beta'' (y+1)^\alpha} \,\mathrm{d}y = \int_{M-1}^\infty
x^{\beta'' y^\alpha} \,\mathrm{d}y = \frac{(-\log(x))^{ -1/\alpha}}{\alpha\beta''^{1/\alpha} }
\Gamma\biggl( \frac{1}{\alpha} , -\beta'' \log(x) (M-1)^\alpha \biggr).
\end{equation}
Using the following asymptotic expansion of $\Gamma(a,b)$ as $b
\rightarrow0$,
%
\begin{equation}
\Gamma(a,b) \sim\Gamma(a) - \frac{b^a}{a}, \qquad a > 0, \mbox{ as }
b \rightarrow0
\end{equation}
(see~\cite{oldham:2009}, formula 45:9:7), and $\log(x) \sim x - 1$ as
$x\rightarrow1$, (\ref{e:rightside}) and (\ref{e:leftside}) are
asymptotic to
%
\begin{equation}
\frac{(1-x)^{-1/\alpha}}{\alpha\beta'^{1/\alpha}}
\Gamma\biggl(\frac{1}{\alpha} \biggr)\quad \mbox{and}\quad \frac{(1-x)^{-1/\alpha}}{\alpha\beta''^{1/\alpha}}
\Gamma\biggl(\frac{1}{\alpha} \biggr),
\end{equation}
respectively. Thus, (\ref{e:inqstring}) implies
%
\begin{equation}\label{e:newstring}
\frac{1}{\alpha\beta''^{1/\alpha}}
\Gamma\biggl(\frac{1}{\alpha}\biggr ) (1-x)^{-1/\alpha} \lesssim
G_{\mathbf{c}}^{(M)} (x) \lesssim\frac{1}{\alpha\beta'^{1/\alpha}}
\Gamma\biggl(\frac{1}{\alpha} \biggr) (1-x)^{-1/\alpha}
\end{equation}
as $x \rightarrow1$. Since everything is tending to $\infty$,
(\ref{e:newstring}) also holds with
$G_{\mathbf{c}}^{(M)}(x)$ replaced with $G_{\mathbf{c}}(x)$. And
finally, since $\varepsilon>0$ was arbitrary, we can let $\varepsilon
\rightarrow0$, making $\beta',\beta'' \rightarrow\beta$, which
implies~(\ref{e:asyGinf}).
\end{pf}

We are now ready to give the L\'{e}vy Khintchine representation of the
Rosenblatt distribution.
\begin{them}\label{t:LKform}
Let $Z_D$ have a Rosenblatt distribution with $0 < D < 1/2$ and let
\[
\mathbf{\lambda}(D)^{-1} = (\lambda_1(D)^{-1},\lambda
_2(D)^{-1},\ldots)
\]
be the sequence of the inverses of the eigenvalues associated to the
integral operator $\sigma(D) \mathcal{K}_D$ defined in
(\ref{e:Koperator}). Then the characteristic function of $Z_D$ can
be written as
%
\begin{equation}\label{e:charfunc}
\phi(\theta) = \mathbb{E}\mathrm{e}^{\mathrm{i} \theta Z_D} = \exp\biggl( \int_0^\infty
(\mathrm{e}^{\mathrm{i} \theta u} - 1 - \mathrm{i} \theta u) \nu_D(u) \,\mathrm{d}u \biggr),
\end{equation}
where $\nu_D$ is supported on $(0,\infty)$ and is given by
%
\begin{equation}\label{e:defofnuLM}
\nu_D(u) = \frac{1}{2u} G_{\mathbf{\lambda}(D)^{-1}} (\mathrm{e}^{-u/2}) =
\frac{1}{2u} \sum_{n=1}^\infty\exp\biggl( -\frac{u}{2 \lambda_n} \biggr),\qquad
u > 0.
\end{equation}
Moreover, $\nu_D$ has the following asymptotic forms as $u \rightarrow
0^+$ and $u \rightarrow\infty$,
%
\begin{eqnarray}
\label{e:asyofnu0}\nu_D(u) &\sim&\frac{ 2^{D/(1-D)} C(D)^{1/(1 - D)}
}{(1-D) }
\Gamma\biggl(\frac{1}{1-D} \biggr) u^{(D-2)/(1-D) }, \qquad u
\rightarrow0,\\
\label{e:asyofnuinf}\nu_D(u) &\sim&\frac{\mathrm{e}^{-u/(2 \lambda_1)}}{2 u}, \qquad u
\rightarrow\infty,
\end{eqnarray}
where $C(D)$ is defined in (\ref{e:defofC}).
\end{them}
\begin{pf}
Let
\[
Z_D^{(M)} = \sum_{n=1}^M \lambda_n (\varepsilon_i^2 - 1).
\]
We have
$Z_D^{(M)} \stackrel{d}{\rightarrow} Z_D$, and since $Z_D^{(M)}$ is a
sum of shifted i.i.d. chi-squared distributions, we can use the L\'
{e}vy Khintchine
representation of a chi-square (\cite{applebaum:2004}, Example
1.3.22), which is a gamma distribution with shape parameter
1$/$2 and scale parameter 2:
%
\begin{equation}
\mathbb{E}\mathrm{e}^{\mathrm{i} \theta Z_D^{(M)}} = \prod_{n=1}^M \mathrm{e}^{\mathrm{i} \theta
\lambda_i
(\varepsilon_i^2 - 1)}
= \prod_{n=1}^M \exp\biggl( -\mathrm{i} \theta\lambda_i + \int_0
^\infty(\mathrm{e}^{\mathrm{i}
\theta u} - 1) \biggl[ \frac{\mathrm{e}^{-u/(2 \lambda_i) }}{2 u}
\biggr] \,\mathrm{d}u \biggr). \label{e:beforeintoint}
\end{equation}
Using
$(1/2)
\int_0^\infty \mathrm{e}^{-u/(2 \lambda)} \,\mathrm{d}u = \lambda,
$
(\ref{e:beforeintoint}) can be rewritten as
\begin{eqnarray*}
& &\prod_{n=1}^M \exp\biggl( \int_0^\infty(\mathrm{e}^{\mathrm{i}
\theta u} - 1 - \mathrm{i} \theta u) \biggl[ \frac{\mathrm{e}^{-u/(2 \lambda_i) }}{2 u}
\biggr] \,\mathrm{d}u \biggr) \\
& & \quad= \exp\biggl( \int_0^\infty(\mathrm{e}^{\mathrm{i}
\theta u} - 1 - \mathrm{i} \theta u) \biggl[ \frac{1}{2u}
G_{\lambda(D)^{-1}}^{(M)}(\mathrm{e}^{-u/2}) \biggr] \,\mathrm{d}u \biggr),
\end{eqnarray*}
where
\[
G_{{\lambda(D)^{-1}}}^{(M)} (x) = \sum_{n=1}^M
x^{\lambda_i^{-1}}.
\]
Now, we let $M \rightarrow\infty$. In order to
justify passing the limit though the integral, notice that
%
\begin{equation}\label{e:boundcf}
\hspace*{-5pt}\biggl| (\mathrm{e}^{\mathrm{i}
\theta u} - 1 - \mathrm{i} \theta u) \biggl[ \frac{1}{2u}
G_{\lambda(D)^{-1}}^{(M)}(\mathrm{e}^{-u/2}) \biggr] \biggr| \leq
\frac{\theta^2}{4} u G_{\lambda(D)^{-1}}^{(M)}(\mathrm{e}^{-u/2})
\leq\frac{\theta^2}{4} u G_{\lambda(D)^{-1}}(\mathrm{e}^{-u/2}),
\end{equation}
where we have used the identity $|\mathrm{e}^{\mathrm{i} z} - 1 - z| \leq\frac{z^2}{2}$ for
$z \in\mathbb{R}$. Notice that (\ref{e:boundcf}) is continuous for
$0 < u < \infty$, and
by (\ref{e:asylambdas}) together with Lemma~\ref{l:gengeosum} using
$\alpha= 1-D$ and $\beta= C(D)^{-1}$, we have
%
\begin{equation}\label{e:lminf}
u G_{\lambda(D)^{-1}}(\mathrm{e}^{-u/2}) \sim u \mathrm{e}^{-u/(2 \lambda_1)}
\qquad\mbox{as } u \rightarrow\infty
\end{equation}
and,
%
\begin{equation}\label{e:lm0}
u G_{\lambda(D)^{-1}}(\mathrm{e}^{-u/2}) \sim C' u(1- \mathrm{e}^{-u/2})^{-1/(1-D)}
\sim C'' u^{D/(1-D)}
\qquad\mbox{as } u \rightarrow0
\end{equation}
for some constants $C'$ and $C''$. Since $0 < \frac{D}{1-D} < 1$,
(\ref{e:lminf}) and (\ref{e:lm0}) imply
that (\ref{e:boundcf}) is integrable on $(0,\infty)$, and hence the
dominated convergence theorem applies and
\[
\mathbb{E}\mathrm{e}^{\mathrm{i} \theta Z_D^{(M)}} \rightarrow\mathbb{E}\mathrm{e}^{\mathrm{i} \theta
Z_D} = \exp\biggl( \int_0^\infty
(\mathrm{e}^{\mathrm{i} \theta u} - 1 - \mathrm{i} \theta u) \biggl[ \frac{1}{2 u}
G_{\lambda(D)^{-1}}(\mathrm{e}^{-u/2}) \biggr] \,\mathrm{d}u \biggr)
\]
which verifies (\ref{e:charfunc}).

The final assertions (\ref{e:asyofnu0}) and (\ref{e:asyofnuinf}) also
follow from
(\ref{e:asylambdas}) and Lemma
\ref{l:gengeosum} with $\alpha= 1-D$ and $\beta= C(D)^{-1}$, since
these imply
%
\begin{eqnarray}
\frac{1}{2u} G_{\lambda(D)^{-1}}(\mathrm{e}^{-u/2}) &\sim&
\frac{1}{2u}
\frac{C(D)^{-1/(1-D)}}{(1-D)} \Gamma\biggl( \frac{1}{1-D}
\biggr) \biggl(\frac{u}{2} \biggr)^{-1/(1-D)} \nonumber
\\[-8pt]
\\[-8pt]
&=& \frac{ 2^{D/(1-D)} C(D)^{1/(1 - D)} }{(1-D) }
\Gamma\biggl(\frac{1}{1-D} \biggr) u^{(D-2)/(1-D)},\qquad
u \rightarrow0\nonumber
\end{eqnarray}
and
%
\begin{equation}
\frac{1}{2u} G_{\lambda(D)^{-1}}(\mathrm{e}^{-u/2}) \sim\frac{1}{2 u}
\mathrm{e}^{-u/(2 \lambda_1)},\qquad u \rightarrow\infty.
\end{equation}
This concludes the proof.\vadjust{\goodbreak}
\end{pf}
\begin{rem*}
 Notice that for any $0 < D < 1/2$, the
L\'{e}vy measure is normalized in the sense that
%
\begin{equation}\label{e:int1}
\int_0^\infty u^2 \nu_D(u) \,\mathrm{d}u = \mathbb{E} Z_D^2 = 1.
\end{equation}
As $D \rightarrow1/2$, $\sigma(D) \rightarrow0$ by
(\ref{e:defofsigma}) $\Rightarrow C(D) \rightarrow0$ by (\ref{e:defofC})
$\Rightarrow\lambda_n(D) \rightarrow0$ by (\ref{e:asylambdas}) $
\Rightarrow\nu_D(u)
\rightarrow0$ by (\ref{e:defofnuLM}) and hence $u^2 \nu_D(u)
\rightarrow0$ for all $u$ positive, but in light of (\ref{e:int1}),
the function $u^2 \nu_D(u)$ approaches a dirac mass at $u=0$. Thus,
as $D \rightarrow1/2$, since $(\mathrm{e}^{\mathrm{i} u \theta} - 1 - \mathrm{i} u \theta)
\rightarrow-1/2 \theta^2$, one gets
\[
\phi(\theta) = \exp\biggl(\int_0^\infty \biggl(\frac{\mathrm{e}^{\mathrm{i} u \theta} - 1
- \mathrm{i} \theta u}{u^2} \biggr) u^2 \nu_D(u) \,\mathrm{d}u\biggr ) \rightarrow
\exp\biggl( -\frac{1}{2} \theta^2 \biggr),
\]
which verifies that $Z_D \stackrel{d}{\rightarrow} N(0,1)$.
\end{rem*}

Understanding the L\'{e}vy measure of a distribution has some
immediate implications pertaining to its probability density function
and distribution function. We will state three such results as
corollaries. The first
is not surprising given that $Z_D$ is an infinite sum of
chi-squared distributions, however it is now easy to prove given
what we know about the L\'{e}vy measure:
\begin{cor}
For $0 < D < 1/2$, the probability density function of $Z_D$ is infinitely
differentiable with all derivatives tending to $0$ as $|x| \rightarrow
\infty$.
\end{cor}
\begin{pf}
Using Proposition 23.8 in~\cite{sato:1999}, this follows as long as
there exists an $\alpha \in(0,2)$ such that
%
\begin{equation}\label{e:smoothcond}
\liminf_{r \rightarrow0} \frac{\int_{[-r,r]} u^2
\nu_D(u) \,\mathrm{d}u}{r^{2-\alpha} }> 0.
\end{equation}
And indeed, (\ref{e:asyofnu0}) implies that for some constant $\tilde{C}(D)$,
%
\begin{equation}
\int_{[-r,r]} u^2 \nu_D(u) \,\mathrm{d}u \sim\tilde{C}(D) \int_0^r
u^{-D/(1-D)} \,\mathrm{d}u = \tilde{C}(D)\biggl (\frac{1 - D}{1-2D} \biggr) r^{2 -
1/(1-D)}.
\end{equation}
Thus, choosing $\alpha= 1/(1-D)$ verifies the result.
\end{pf}

The second corollary gives a simple bound on the left-hand tail of the
CDF of
$Z_D$. The proof is similar to Proposition 9.5(ii) in
\cite{steutel:2004}.
\begin{cor}[(Left tail)]
Let $Z_D$ denote the Rosenblatt distribution. Then,
%
\begin{equation}
\mathrm{P}[Z_D < -x] \leq\exp\bigl( -\tfrac{1}{2} x^2 \bigr),\qquad x
> 0.
\end{equation}
\end{cor}
\begin{pf}
For any $M \geq1$, recall the random variable $Z^{(M)}_D$ defined in the
proof of Theorem~\ref{t:LKform}. Applying Markov's inequality we see
for any $x > 0$ and $s > 0$,
%
\begin{eqnarray}\label{e:beforeBoundCOR2}
\mathrm{P}\bigl( Z^{(M)}_D \leq- x\bigr) &\leq& \mathrm{e}^{-s x} \mathbb{E} \mathrm{e}^{- s
Z^{(M)}_D} \nonumber
\\[-8pt]
\\[-8pt]
&=& \mathrm{e}^{-s x} \exp\biggl( \int_0^\infty(\mathrm{e}^{- s u} - 1 +s u)
\,\mathrm{d}\nu^{(M)}(u) \,\mathrm{d}u \biggr),\nonumber
\end{eqnarray}
where the last equality follows since $Z^{(M)}_D$ is a finite sum of
weighted chi-square distributions, and hence has a moment generating
function defined for $s<s_0$ for some $s_0 > 0$. Using the inequality
$\mathrm{e}^{- s u} - 1 +s u \leq\frac{1}{2} s^2 u^2$, and since
$\nu^{(M)}(u)$ increases in $M$,
(\ref{e:beforeBoundCOR2}) becomes
%
\begin{eqnarray}\label{e:boundofMGF}
 \mathrm{e}^{-s x} \exp\biggl( \int_0^\infty(\mathrm{e}^{- s u} - 1 +s u)
\,\mathrm{d}\nu^{(M)}(u) \,\mathrm{d}u \biggr) &\leq& \mathrm{e}^{-s x} \exp\biggl(\frac{s^2}{2}
\int_0^\infty u^2 \nu^{(M)}(u) \,\mathrm{d}u \biggr) \nonumber
\\
& \leq& \mathrm{e}^{-s x} \exp\biggl(\frac{s^2}{2}
\int_0^\infty u^2 \nu(u) \,\mathrm{d}u \biggr)\\
&\leq &\exp\biggl( -s x + \frac{1}{2} s^2 \biggr)\nonumber
\end{eqnarray}
since $\operatorname{Var} Z_D = \int_0^\infty u^2 \nu(u) \,\mathrm{d}u =1$. The
minimum of (\ref{e:boundofMGF}) over $s$ is attained at $s = x$. Thus,
%
\begin{equation}
\mathrm{P}\bigl( Z^{(M)}_D \leq- x\bigr) \leq\exp\bigl(-\tfrac{1}{2} x^2 \bigr)
\end{equation}
which finishes the proof.
%
%
\end{pf}

We obtain a result similar to Theorem 2 in
\cite{albin:1998} involving the rate of decay of the distribution function.
\begin{cor}[(Right tail)]
Let $Z_D$ denote the Rosenblatt distribution. Then for $\alpha> 0$,
\[
\lim_{u \rightarrow\infty} \frac{\mathrm{P}[Z_D > u +
\alpha]}{\mathrm{P}[Z_D > u]} =
\mathrm{e}^{-\alpha/(2 \lambda_1)},
\]
where $\lambda_1$ is the largest eigenvalue of $\sigma(D) \mathcal{K}_D$
defined in (\ref{e:Koperator}).
\end{cor}
\begin{pf}
For $u \geq1$, define $\bar{A}_{\nu_D}$ as
\[
\bar{A}_{\nu_D}(u) = \frac{\int_u^\infty\nu_D(v) \,\mathrm{d}v }{ \int
_1^\infty
\nu_D(v) \,\mathrm{d}v}.
\]
By Theorem 1 in~\cite{braverman:2005}, it suffices to show
%
\begin{equation}\label{e:showme1}
\lim_{u \rightarrow\infty} \frac{\bar{A}_{\nu_D}(u + \alpha)}{
\bar{A}_{\nu_D}(u)} =
\mathrm{e}^{-\alpha/(2 \lambda_1)},
\end{equation}
and
%
\begin{equation}\label{e:showme2}
\int_0^\infty \mathrm{e}^{u/(2 \lambda_1)} \bar{A}_{\nu_D}(\mathrm{d}u) =
\infty.
\end{equation}
Both follow from (\ref{e:asyofnuinf}) in Theorem~\ref{t:LKform}.
Indeed, as $u \rightarrow\infty$,
%
\begin{equation}\label{e:asyofUIGF}
\int_u^\infty\nu_D(v) \,\mathrm{d}v \sim\int_u^\infty\frac{1}{2 v}
\mathrm{e}^{-v/(2 \lambda_1)} \,\mathrm{d}v \sim\frac{\lambda_1}{u}
\mathrm{e}^{-u/(2 \lambda_1)}
\end{equation}
which follows because
\[
-\frac{\mathrm{d}}{\mathrm{d}x} \frac{\lambda_1}{u}
\mathrm{e}^{-u/(2 \lambda_1)} = \frac{1}{2 u} \mathrm{e}^{-u/(2 \lambda_1)} +
\frac{\lambda_1}{u^2} \mathrm{e}^{-u/(2 \lambda_1)} \sim\frac{1}{2 u} \mathrm{e}^{-u/(2
\lambda_1)},
\]
and thus integrating both sides implies (\ref{e:asyofUIGF}).
From this, (\ref{e:showme1}) follows. Also,
\begin{eqnarray*}
\int_0^\infty \mathrm{e}^{u/(2 \lambda_1)} \bar{A}_{\nu_D}(\mathrm{d}u) &=&
\frac{1}{ \int_1^\infty
\nu_D(v) \,\mathrm{d}v} \int_0^\infty \mathrm{e}^{u/(2 \lambda_1)} \Biggl(
\frac{1}{2u} \sum_{n=1}^\infty \mathrm{e}^{-u/(2 \lambda_n)} \Biggr) \,\mathrm{d}u\\
&=& \frac{1}{ \int_1^\infty
\nu_D(v) \,\mathrm{d}v} \int_0^\infty\frac{1}{2u} \Biggl( 1 + \sum_{n=1}^\infty
\mathrm{e}^{-u/(2
\lambda_n) + u/(2 \lambda_1)} \Biggr) \,\mathrm{d}u = \infty,
\end{eqnarray*}
since the integrand on the right is asymptotic to $1/(2u)$. This verifies
(\ref{e:showme2}).
\end{pf}
%
\section{Approximating the distribution of the tail of the series}
The representation of the Rosenblatt distribution $Z_D$ as an infinite sum
of shifted chi-squared distributions (\ref{e:chisum}) has proven to be
quite useful for
obtaining theoretical results about the distribution. In this
section, we aim to take advantage of this representation to compute
the CDF and PDF of the Rosenblatt distribution.

For $M \geq1$, define $X_M$
and $Y_M$ as
%
\begin{equation}\label{e:decomp}
Z_D = \sum_{n=1}^{M-1} \lambda_n(\varepsilon_n^2 - 1) + \sum
_{n=M}^\infty
\lambda_n(\varepsilon_n^2 - 1) := X_M + Y_M.
\end{equation}
Notice that $Y_M$ has mean $0$ and variance
%
\begin{equation}\label{e:defofsigM7}
\sigma_M^2 := \mathbb{E} Y_M^2 = 2 \sum_{n=M}^\infty\lambda_n^2.
\end{equation}
Notice that from Theorem~\ref{t:evals}, as $ M\rightarrow\infty$,
%
\begin{equation}\label{e:siglargeM}
\sigma_M^2 \sim2C(D)^2 \sum_{n=M}^\infty n^{2D - 2} \sim 2 C(D)^2
(1-2D)^{-1} M^{2D - 1},
\end{equation}
where we approximated the sum $\sum_{n=M}^\infty\lambda_n^2$ with
an integral. This suggests that $X_M$ alone
is not a very good approximation of $Z_D$ since this variance tends to
$0$ slowly with~$M$, especially when $D$ is close to $1/2$. As an
alternative, we will see below that $Y_M$ is approximated by a normal
distribution as $M \rightarrow\infty$. By taking advantage of this
property, we can obtain accurate approximations of the distribution of $Z_D$.

In~\cite{veilletteEW:2010}, random variables expressed as infinite
sums of gamma and hence chi-squared distributions were studied. It was
shown in particular that the asymptotic form of the $\lambda_n$'s in
the expansion
implies that the distribution of the tail $Y_M$ is approximately normal
as $M \rightarrow\infty$.

Let $\kappa_{k,M}$ be the normalized cumulants of $Y_M$. These are
given by
%
\begin{equation}\label{e:defofkapMR}
\kappa_{k,M} = 2^{k-1} (k-1)! \sigma_M^{-k} \sum_{n=M}^\infty
\lambda_n^k.
\end{equation}
Notice that from the asymptotics (\ref{e:asylambdas}) of $\lambda_n$,
\begin{eqnarray*}
\kappa_{k,M} &\sim& 2^{k-1} (k-1)! \cdot 2^{-k/2} (1 - 2D)^{k/2}
M^{ k/2 - kD} \cdot\int_M^\infty
n^{k(D-1)} \,\mathrm{d}n \\
&=& \frac{(k-1)!}{2} \frac{(2 - 4D)^{k/2}}{k-kD-1} M^{1- k/2},
\end{eqnarray*}
which is indeed equal to 1 when $k=2$.

The convergence to normality of $Y_M$ is implied by the following
Berry--Essen bound proved in~\cite{veilletteEW:2010}.
\begin{them}\label{t:berryessen}
Let $Y_M$ and $\sigma_M$ be defined as in (\ref{e:decomp}) and
(\ref{e:defofsigM7}). Then\footnote{Multiplicative factors: they
are denoted $\sigma(D)$ in (\ref{e:scalerv}) and $\sigma_M^{-1}$ in
(\ref{e:berryessen7}). Using (\ref{e:ROSseries})
and (\ref{e:defofsigM7}), compare also (\ref{e:cumulantformula}) with
(\ref{e:defofkapMR}).}\tsup{,}\footnote{The constant 0.7056 appearing in
this inequality is the smallest known to date, see~\cite{shevtsova:2006}.},
%
\begin{equation}\label{e:berryessen7}
\sup_{x \in\mathbb{R}} |\mathrm{P}[\sigma_M^{-1} Y_M \leq x] -
\Phi(x)|
\leq 0.7056\kappa_{3,M} = \mathrm{O}(M^{-1/2}).
\end{equation}
\end{them}

While this result gives an idea of the distribution of $Y_M$ as $M
\rightarrow\infty$, the bound on the right-hand side of
(\ref{e:berryessen7}) may not be satisfactory for small $M$. In order
to better approximate the CDF of $Y_M$, we will use an Edgeworth
expansion which is also considered in~\cite{veilletteEW:2010}. These
expansions use the higher cumulants to better approximate the CDF of
$Y_M$ and give a faster convergence than what we see in (\ref{e:berryessen7}).

Recall the \textit{Hermite polynomials}, which are defined as
\[
H_k(x) = (-1)^k \mathrm{e}^{x^2/2} \frac{\mathrm{d}^k}{\mathrm{d}x^k} \mathrm{e}^{-x^2/2},\qquad k \geq0.
\]
The first few are $H_2(x) = x^2 -1$, $H_3(x) =
x^3 - 3x$, $H_4(x) = x^4 - 6x^2 + 3$ and $H_5(x) = x^5 - 10x^3 + 15x$.
A simple induction gives
%
\begin{equation}\label{e:diffHermite}
\frac{\mathrm{d}}{\mathrm{d}x} H_k(x) \phi(x) = -H_{k+1}(x) \phi(x),
\end{equation}
where $\phi$ is the standard normal PDF. The following theorem is also
proved in~\cite{veilletteEW:2010}.
\begin{them}\label{t:RosEWexpansion}
The CDF of the tail $Y_M$ satisfies
%
\begin{eqnarray}\label{e:EWexpanROS}
\mathrm{P}[\sigma_M^{-1} Y_M \leq x] &=&\Phi(x) - \phi(x) \Biggl\{
\sum_{\eta(N)} \Biggl[ \prod_{m=1}^N \frac{1}{k_m!} \biggl(
\frac{\kappa_{m,M}}{m!} \biggr)^{k_m} \Biggr]
H_{\zeta(k_3,\ldots,k_N)} (x) \Biggr\}\nonumber
\\[-8pt]
\\[-8pt]
 &&{}+ \mathrm{O} \bigl( M^{-(N-1)/2} \bigr),\nonumber
\end{eqnarray}
where $\kappa_{k,M}$ is defined in (\ref{e:defofkapMR}), $\eta(N)$
denotes all $k_3,k_4,\ldots,k_N$ such that
%
\begin{equation}
1 \leq k_3 + 2 k_4 + \cdots+ (N-2) k_n \leq N-2
\end{equation}
and
%
\begin{equation}
\zeta(k_3,\ldots,k_n) = 3 k_3 + 4 k_4 + \cdots+ N k_N - 1.
\end{equation}
\end{them}

Notice by taking a derivative of (\ref{e:EWexpanROS}) and using (\ref
{e:diffHermite}), the PDF of $Y_M$ is approximated by
%
\begin{equation}\label{e:RosPDFapprox}
f_{\sigma_M^{-1} Y_M}(x) = \phi(x) \Biggl[ 1 + \sum_{\eta(N)} \Biggl[ \prod
_{m=1}^N \frac{1}{k_m!} \biggl(
\frac{\kappa_{m,M}}{m!} \biggr)^{k_m} \Biggr]
H_{\zeta(k_3,\ldots,k_N) + 1} (x) \Biggr] + \mathrm{O} \bigl(
M^{-(N-1)/2} \bigr).\hspace*{4pt}
\end{equation}

These expansions allows us to compute the CDF and PDF of $Y_M$ (and
hence of
$Z_D$) to high accuracy. We provide the details and results of this
method in Section~\ref{s:computeCDF}.
%
\section{Obtaining the eigenvalues, cumulants and moments
numerically}\label{s:compute}
%

There exists an extensive literature regarding the problem of
approximating eigenvalues of
integral operators like $\mathcal{K}_D$, see for instance~\cite{spence:1975,atkinson:1967,sloan:1976,ahues:2001}, or~\cite{chen:2009}.
Many, if not all, of these methods boil down to approximating $\mathcal
{ K}_D$ with a finite-dimensional linear operator, a technique often
referred to as the
``Nystr{\"o}m method'' (see~\cite{spence:1978} and~\cite{atkinson:han:2009}).
We shall approximate $\mathcal{ K}_D$ by a $(J+1) \times(J+1)$ matrix
$\mathbf{K}_D$, defined as follows.
%

Fix $J > 0$ and choose nodes $0 = x_0
< x_1 < \cdots< x_J = 1$. For such $f$, let $\mathbf{f} =
(f(x_0),\ldots, f(x_{J}))$
and define
%
\begin{equation}
f_J(x) = f(x_{j-1})\frac{x- x_j}{x_{j-1} - x_j} + f(x_j) \frac{x -
x_{j-1}}{x_j - x_{j-1}}\qquad \mbox{for } x \in[x_{j-1},x_j].
\end{equation}
Then, for any $i=0,1,\ldots,J$,
\[
(\mathcal{K}_D f_J)(x_i) =
\sum_{j=1}^J \int_{x_{j-1}}^{x_j} \biggl( f(x_{j-1})\frac{u- x_j}{x_{j-1}
- x_j} + f(x_j) \frac{u -
x_{j-1}}{x_j - x_{j-1}} \biggr) |x_i - u|^{-D} \,\mathrm{d}u = (\mathbf{K}_D \mathbf{f})_i,
\]
where
$J$ indicates the level of approximation.
For more details on the matrix $\mathbf{K}_D$, see Sections~A and B of
\supp.
By
taking the eigenvalues of the matrix~$\mathbf{K}_D$, we can
approximate the $M$ largest eigenvalues of $\mathcal{K}_D$ where $M
\ll
J$.

%
\begin{table}
\tablewidth=\textwidth
\tabcolsep=0pt
  \caption{First 10 eigenvalues of $\mathcal{K}_D$ for various $D$. Note
that for $D=0$, we have $\lambda_1 = 1$ and $\lambda_n = 0$ for $n
\geq2$.
As $D \rightarrow1/2$, $\lambda_n \rightarrow0$}\label{tab:evals}
\begin{tabular*}{\textwidth}{@{\extracolsep{\fill}}lllll@{}}
\hline
$n$ & $D = 0.1$ & $D = 0.2$ & $D = 0.3$ & $D = 0.4$ \\
\hline
\phantom{0}1 & 0.70200 & 0.68130 & 0.63050 & 0.51250 \\
\phantom{0}2 & 0.05648 & 0.11260 & 0.16040 & 0.17840 \\
\phantom{0}3 & 0.03477 & 0.07341 & 0.11070 & 0.13010 \\
\phantom{0}4 & 0.02453 & 0.05411 & 0.08512 & 0.10420 \\
\phantom{0}5 & 0.01947 & 0.04409 & 0.07119 & 0.08948 \\
\phantom{0}6 & 0.01600 & 0.03710 & 0.06129 & 0.07879 \\
\phantom{0}7 & 0.01374 & 0.03240 & 0.05446 & 0.07121 \\
\phantom{0}8 & 0.01198 & 0.02872 & 0.04903 & 0.06512 \\
\phantom{0}9 & 0.01069 & 0.02596 & 0.04489 & 0.06037 \\
10 & 0.009629 & 0.02366 & 0.04140 & 0.05635 \\
\hline
\end{tabular*}
\end{table}

In order to test that the approximate eigenvalues we are obtaining are
accurate, we have three methods:
\begin{enumerate}[3.]
\item[1.] Check for numerical convergence, that is, for $J$ large, does
increasing $J$ lead to a negligible change in the approximations of
$\lambda_n$. Thus, to approximate $\lambda_n$, we increased $J$
until sufficient convergence was met. By ``sufficient
convergence,'' we mean that $J$ is increased by multiples of $200$
until the
values of the $\lambda_n$'s no longer changed in the 4
significant decimal digit. Table~\ref{tab:evals} gives the results
of this method applied to the first 10 values of $\lambda_n$'s. For
more values, see \supp.

\item[2.] Compare $\lambda_n$ with the asymptotic formula given by
Theorem~\ref{t:evals} for $n$ large. The first 30 $\lambda_n$ are
approximated and plotted on a log scale in \supp and
compared to the asymptotic formula in Theorem~\ref{t:evals}. For
large $n$,
our approximations appear to be in agreement with the asymptotic
formula as $n
\rightarrow\infty$. In fact, the asymptotic
formula for the $\lambda_n$ is a good approximation even for moderate
values of $n$ (about $n \geq20$ for $D=0.1,0.2,0.3,0.4$).

\item[3.] The cumulants $\kappa_{k}$ of $Z_D$ are given in terms of the
eigenvalues as
%
\begin{equation}\label{e:kappakagain}
\kappa_k = 2^{k-1} (k-1)! \sum_{n=1}^\infty\lambda_n^k = 2^{k-1}
(k-1)! \Biggl[ \sum_{n=1}^{M-1} \lambda_n^{k} + \sum_{n=M}^\infty
\lambda_n^k \Biggr].
\end{equation}
Since $\kappa_2 = 1$, and $\kappa_3, \kappa_4$ can be
computed numerically exactly, we can
compute the absolute error made by approximating the cumulants
$\kappa_2,\kappa_3,\kappa_4$ by using (\ref{e:kappakagain}). In order
to compute the sums on the right-hand side of (\ref{e:kappakagain}),
we used the first $M=50$ values of $\lambda_n$ using the matrix
$\mathbf{K}_D$, and for $n > M$,
the asymptotic formula given in Theorem~\ref{t:evals}. The results of
this test are given in a table in \supp. Because we are
using the asymptotic formula for $n$ large, we expect some error, but
as it turns out, the absolute errors are very small compared to
the size of $\kappa_k$, $k=2,3,4$.
\end{enumerate}

One can obtain also higher cumulants and moments of $Z_D$.
The cumulants of $Z_D$ are expressed in terms of the functions
$G_{k,D}$ by~\ref{e:cumulantformula} and
Proposition~\ref{p:computeint}. These functions are defined
recursively in (\ref{e:defofGkD}). By using
$\mathbf{K}_D$ to approximate the functions $G_{k,D}$, as explained in
Section C of \supp, one obtains
approximations for the cumulants. One derives the corresponding moments
by applying
(\ref{e:momentsBell}). We have tabulated the first 8
moments and cumulants of the Rosenblatt distribution for various $D$
in \supp.
%
\section{Obtaining the CDF numerically}\label{s:computeCDF}
We now present a technique for computing the CDF and PDF of the
Rosenblatt distribution. Computing the distribution of $Z_D=X_M+ Y_M$
(see (\ref{e:decomp})) requires
three steps:
\begin{enumerate}[3.]
\item[1.] Approximate the eigenvalues $\lambda_n$ for
$n=1,2,\ldots,M$ for some $M > 0$.

\item[2.] Compute separately the CDF of $X_M$ and the PDF of $Y_M$ in the
decomposition (\ref{e:decomp}).

For $X_M$, methods exist already for
accurately computing the CDF of a finite sum of
chi-squared distributions, see for instance methods based on Laplace transform
inversion,~\cite{veillette:2010,castano:2005}, or Fourier
transform inversion,~\cite{abate:1992}.

For $Y_M$, use an Edgeworth
expansion of order $N \geq2$ which is found in Theorem
\ref{t:RosEWexpansion}. To compute $\sigma_M$ in (\ref
{e:defofsigM7}), take
$
\sigma_M = (1 - 2 \sum_{n=1}^{M-1} \lambda_n^2)^{1/2}.
$
The Edgeworth expansion in Theorem~\ref{t:RosEWexpansion} also
involves the
cumulants $\kappa_{k,M}$, $k \geq1$ of $Y_M$. If $M$ is
sufficiently large, $\kappa_{k,M}$ can be approximated using
the asymptotic formula given in Theorem~\ref{t:evals}:
%
\begin{eqnarray}\label{e:useHurwitz}
\kappa_{k,M} &=& \sigma_M^{-k} \Biggl( 2^{k-1} (k-1)!
\sum_{n=M}^{\infty} \lambda_n^k\Biggr )\approx \sigma_M^{-k} \Biggl( 2^{k-1} (k-1)!
\sum_{n=M}^{\infty} C(D)^k n^{k(D-1)} \Biggr)\hspace*{5pt}\nonumber\quad
\\[-8pt]
\\[-8pt]
&=& 2^{k-1} (k-1)! \sigma_M^{-k} C(D)^k \zeta\bigl(k(1-D),M\bigr),\nonumber\quad
\end{eqnarray}
where $\zeta(s,M) = \sum_{n=M}^\infty n^{-s}$ denotes the \textit
{Hurwitz Zeta function}
(\cite{oldham:2009}, Chapter 64). Using (\ref{e:useHurwitz})
introduces a small amount of error since we are approximating
$\lambda_n$ for $n$ large, however we will see below this is
negligible for $M$ large.

\item[3.] Finally, the CDF of $Z_D$ is given by a convolution
\[
F_{Z_D}(x) = \int_{-\infty}^\infty F_{X_M}(x - y) f_{Y_M}(y) \,\mathrm{d}y,
\]
where $F_{X_M}$ is the CDF of $X_M$ and $f_{Y_M}$ is the PDF of
$Y_M$. We compute this convolution in MATLAB using standard numerical
integration techniques.
\end{enumerate}
%
%
\begin{table}
\tablewidth=\textwidth
\tabcolsep=0pt
  \caption{Approximation of $\mathrm{P}[Z_{0.3} \leq0]$ for various $M$
and $N$ values. This shows that one gets an excellent approximation
with relatively small values of $M$ and $N$}\label{tab:incMN}

\begin{tabular*}{\textwidth}{@{\extracolsep{\fill}}lllll@{}}
\hline
$N$ & $M=10$ & $M = 20$ & $ M = 30$ & $M = 50$ \\
\hline
2 & 0.616883& 0.616909 & 0.616909 & 0.616907 \\
3 & 0.616817 & 0.616878 & 0.616890 & 0.616896 \\
4 & 0.616885 & 0.616898 & 0.616899 & 0.616899 \\
5 & 0.616895 & 0.616900 & 0.616900 & 0.616900 \\
6 & 0.616895 & 0.616900 & 0.616900 & 0.616900 \\
\hline
\end{tabular*}
\end{table}

The choice of $M$ (number of chi-squared distributions) and $N$
(order of Edgeworth expansion) was determined by increasing both until
the value of the CDF changed by less than $10^{-5}$. In
Table~\ref{tab:incMN}, we show the effects of increasing $M$ and $N$ for the
case of $D = 0.3$ and $x = 0$. Observe that for fixed
$M$, the values of the CDF converge rapidly as $N$ increases.
Nevertheless, if $M$ is small $(M \leq10)$ the approximation will have
a slight error since we are approximating the $\kappa_{k,M}$ using
(\ref{e:useHurwitz}). Moreover, by choosing $N=2$ and $M \geq10$, the
only noticeable improvement by increasing $N$ happens in the fifth decimal
place and beyond.
%
\begin{table}[b]
\tablewidth=\textwidth
\tabcolsep=0pt
\caption{Various quantiles of the Rosenblatt distribution for
selected values of $D$}\label{tab:quantiles}
\begin{tabular*}{\textwidth}{@{\extracolsep{\fill}}ld{2.4}d{2.4}d{2.4}d{2.4}d{2.4}@{}}
\hline
Quantile& \multicolumn{1}{l}{$D=0.1$} & \multicolumn{1}{l}{$D=0.2$}
& \multicolumn{1}{l}{$ D=0.3 $}& \multicolumn{1}{l}{$ D=0.4 $}
& \multicolumn{1}{l@{}}{$ D=0.45$} \\
\hline
$0.01 $& -0.8472 & -1.0567 & -1.3546 & -1.7838 & -2.0603\\
$0.025 $& -0.8142 & -0.9827 & -1.7669 & -1.5536 & -1.7639 \\
$0.05 $& -0.7808 & -0.9122 & -1.0958 & -1.3493 & -1.5051 \\
$0.10 $& -0.7340 & -0.8201 & -0.9419 & -1.1053 & -1.6462 \\
$0.25 $& -0.6200 & -0.6277 & -0.6479 & -0.6713 & -0.6789 \\
$0.50 $& -0.3622 & -0.3055 & -0.2329 & -0.1332 & -0.0673 \\
$0.75 $& 0.2370 & 0.2781 & 0.3666 & 0.5059 &  0.5955 \\
$0.90 $& 1.2047 & 1.2031 & 1.2110 & 1.2417 &  1.2673 \\
$0.95 $& 2.0015 &  1.9726 & 1.9114 & 1.8022 &  1.7262 \\
$0.975 $& 2.8312 &  2.7759 & 2.6483 & 2.3858 &  2.1774 \\
$0.99 $& 3.9618 &  3.8718 & 3.6579 & 3.1909 &  2.7892 \\
\hline
\end{tabular*}
\end{table}

\begin{figure}

\includegraphics{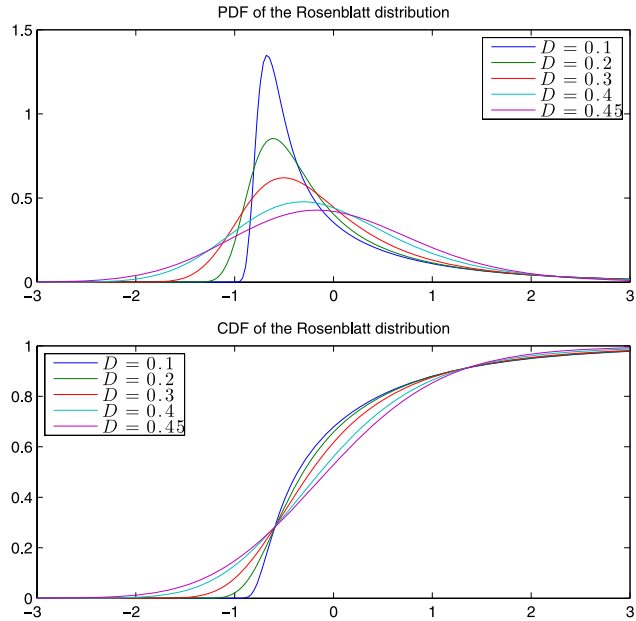}

  \caption{Plots of the PDF and CDF of $Z_D$ for various $D$. The CDF
with the steepest slope
and the PDF with the highest mode correspond to $D=0.1$.}\label{f:PDFCDFZ}
\end{figure}

We have tabulated in Table~\ref{tab:quantiles} quantiles of $Z_D$ for
various $D$.
These are useful for obtaining confidence intervals. To obtain these
values, we solved the equation $F_{Z_D}(x) = q$ for various quantiles
$q$ in MATLAB. To compute the CDF $F_{Z_D}(x)$, we fixed $N=6$, and
increased $M$ in increments of 10 until the approximation
of the CDF changed by less than $10^{-5}$. For the CDF values, see
Table 9
in \supp. We have also plotted the PDF and CDF in Figure~\ref{f:PDFCDFZ}.
%

\section*{Acknowledgement}
This work was partially supported
by the NSF Grants DMS-07-06786 and DMS-10-07616 at Boston University.

\begin{supplement}
\sname{Supplemental article}
\stitle{Supplement to \emph{Properties and numerical evaluation of
the Rosenblatt distribution}}
\slink[doi]{10.3150/12-BEJ421SUPP}
\sdatatype{.pdf}
\sfilename{BEJ421\_supp.pdf}
\sdescription{The supplement~\cite{veillette-Rosen-supp:2011} to this
article details the approximation of the integral operator $\mathcal{K}_D$
and the computation of the cumulants, moments and CDF of $Z_D$. It also
contains an extensive table of the CDF of $Z_D$ and a guide to the software.}
\end{supplement}


\printhistory

\end{document}